\documentclass[12pt]{article}
\usepackage{amsmath, amssymb}

\usepackage{latexsym}
\usepackage{theorem}
\usepackage{cite}
\theorembodyfont{\upshape}
\newtheorem{Def}{Definition}[section]
\newtheorem{Thm}[Def]{Theorem}
\newtheorem{Lem}[Def]{Lemma}
\newtheorem{Cor}[Def]{Corollary}
\newtheorem{Prop}[Def]{Proposition}
\newtheorem{Rem}[Def]{Remark}
\newtheorem{Ex}[Def]{Example}

\def\R{\mathbb{R}}

\def\cD{{\cal D}}

\def\cL{{\cal L}}

\def\rfe#1{(\ref{#1})}

\makeatletter
\@addtoreset{equation}{section}
\makeatother

\def\commentmark{%
 \renewcommand{\thefootnote}{$\ast$}%
 \footnotemark%
 \renewcommand{\thefootnote}{\arabic{footnote}}}
\def\commenttext#1{%
 \renewcommand{\thefootnote}{$\ast$}%
 \footnotetext{#1}%
 \renewcommand{\thefootnote}{\arabic{footnote}}}

\newcommand{\qed}{\hfill $\Box$ \\}

\begin{document}

\title{
Semismall perturbations,\\ semi-intrinsic ultracontractivity,\\
and integral representations of nonnegative \\
solutions for parabolic equations\commentmark{}} 
\commenttext{%
}

\author{Pedro J. Mendez-Hernandez \\
Escuela de Matem\'atica,
 Universidad de Costa Rica\\
San Jos\'e, Costa Rica\\
e-mail: pedro.mendez@ucr.ac.cr\and
Minoru Murata \\
Department of Mathematics, Tokyo Institute of Technology, \\
Oh-okayama, Meguro-ku, Tokyo, 152-8551 Japan \\
e-mail: minoru3@math.titech.ac.jp \\
}
\date{ }

\maketitle

\begin{abstract}
We consider nonnegative solutions of a parabolic equation in a cylinder
$D \times I$, where $D$ is a noncompact domain 
of a Riemannian manifold and $I =(0,T)$ with $0 < T \le \infty$ or
$I=(-\infty,0)$.
Under the assumption [SSP] (i.e., the constant function $1$ is a semismall 
perturbation of the associated elliptic operator on $D$),
we establish an integral representation theorem of nonnegative solutions:
In the case $I =(0,T)$,
any nonnegative solution is represented uniquely by an integral on
$(D \times \{ 0 \}) \cup (\partial_M D \times [0,T))$, where
$\partial_M D$ is the Martin boundary of $D$ for the elliptic operator;
and in the case $I=(-\infty,0)$, 
any nonnegative solution is represented uniquely by the sum
of an integral on $\partial_M D \times (-\infty,0)$ and a constant
multiple of a particular solution.
We also show that [SSP] implies the condition [SIU] 
(i.e., the associated heat kernel is semi-intrinsically 
ultracontractive). 
\end{abstract}

\section{Introduction} \label{sec:1}

This paper is a continuation of \cite{Mu07}. It is 
concerned with integral representations
of nonnegative solutions to parabolic equations and 
perturbation theory for elliptic operators.

We consider nonnegative solutions of a parabolic equation
\begin{equation} \label{eq:1.1}
(\partial_t + L)u=0 \quad \text{in} \quad D \times I, 
\end{equation}
where $\partial_t = \partial / \partial t$,
$L$ is a second order elliptic operator on
a noncompact domain $D$ of a
Riemannian manifold $M$, and  $I$ is a time interval: 
$I=(0,T)$ with $0 <T \le \infty$ or $I=(-\infty,0)$.

During the last few decades, much attention has been paid to 
the structure of all nonnegative solutions to a parabolic equation,
perturbation theory for elliptic operators, and their relations.
(See \cite{Ai98}, \cite{AM96}, \cite{An97}, \cite{AT92}, \cite{Aro68}, 
\cite{BB92}, \cite{GH98}, \cite{IM01}, \cite{KT85}, 
\cite{LP94}, \cite{MT84}, \cite{Mu93}, \cite{Mu94}, 
\cite{Mu95a}, \cite{Mu95b}, \cite{Mu96}, \cite{Mu97}, \cite{Mu02}, 
\cite{Mu03}, \cite{Mu05}, \cite{Mu07}, 
\cite{Pi88}, \cite{Pi89}, \cite{Pi99}, \cite{Pin95}, 
\cite{Sm96}, \cite{To07}.)
Among others, Murata \cite{Mu07} has established integral representation 
theorems of nonnegative solutions to the equation (\ref{eq:1.1}) under the 
condition [IU] (i.e., intrinsic ultracontractivity) on the minimal 
fundamental solution $p(x,y,t)$ for (\ref{eq:1.1}). Furthermore, he has 
shown that [IU] implies [SP] (i.e., the constant function $1$ is a
small perturbation of $L$ on $D$). It is known (\cite{Mu97}) that
[SP] implies [SSP] (i.e., $1$ is a semismall perturbation of $L$ on $D$).

In this paper, we show that [SSP] implies [SIU]
(i.e., semi-intrinsic ultracontractivity) and give integral representation 
theorems of nonnegative solutions to (\ref{eq:1.1}) under the 
condition [SSP].  We consider that [SSP] is one of the weakest possible condition for 
getting "explicit" integral representation theorems.

Now, in order to state our main results, we fix notations and recall
several notions and facts.
Let $M$ be a connected separable $n$-dimensional
smooth manifold with Riemannian metric of class $C^0$.
Denote by $\nu $ the Riemannian measure on $M $. 
$T_x M$ and $TM$ denote the tangent space
to $M$ at $x \in M$ and the tangent bundle, respectively.
We denote by $\mathrm{End}(T_x M)$ and $\mathrm{End}(TM)$ the set
of endmorphisms in $T_x M$ and the corresponding bundle, respectively.
The inner product on $TM$ is denoted by $\langle X, Y \rangle$,
where $X,Y \in TM$; and $|X| = \langle X, X \rangle^{1/2}$.
The divergence and gradient with respect to the metric on $M$
are denoted by $\mathrm{div}$ and $\nabla$, respectively.
Let $D$ be a noncompact domain of $M$.
Let $L$ be an elliptic differential operator on $D$ of the form
\begin{equation} \label{eq:1.2}
L  u = - m ^{-1} \mathrm{div} (m  A  \nabla u ) + V  u,
\end{equation}
where $m $ is a positive measurable function on $D $ such that
$m$ and $m ^{-1}$ are bounded on any compact subset of $D$,
$A $ is a symmetric measurable section on $D $ of $\mathrm{End}(TM )$,
and $V $ is a real-valued measurable function on $D $ such that
\begin{equation*} 
V  \in L_{\mathrm{loc}}^p (D ,m d\nu )
\quad \mbox{ for some } p > \max (\frac{n}{2}, 1).
\end{equation*}
Here $L_{\mathrm{loc}}^p (D ,m d\nu )$ is the set of 
real-valued functions on $D $ locally $p$-th integrable
with respect to $m d\nu $.
We assume that $L $ is locally uniformly elliptic on $D $, i.e.,
for any compact set $K$ in $D $ there exists a positive constant
$\lambda$ such that
\begin{equation*} 
\lambda |\xi|^2 \le \langle A_x \xi, \xi \rangle \le \lambda^{-1} |\xi|^2,
\quad x \in K, ~(x, \xi) \in TM .
\end{equation*}
We assume that the quadratic form 
$Q$ on $C_0^{\infty } (D)$ defined by 
\begin{equation*}
Q[u] =\int_D 
(\langle A\nabla u,\nabla u \rangle
+ V u^{2}) m d\nu 
\end{equation*}
is bounded from below, and put
\begin{equation*}
\lambda_0 = \inf \left\{ Q[u]; u \in C_0^{\infty } (D), \ \ 
\int_D   u^{2} m d\nu = 1 \right\}.
\end{equation*}
Then, for any $a < \lambda_0$, $\  (L-a,D)$ is subcritical, i.e.,
there exists the (minimal positive) Green function of $L-a$ on $D$.
We denote by
$L_D$ the selfadjoint operator in $L^2(D;md\nu)$ 
associated with the closure of $Q$.
The minimal fundamental solution for (\ref{eq:1.1}) is denoted by
$p(x,y,t)$, which is equal to the integral kernel of the semigroup
$e^{-t L_D}$ on $L^2(D,md\nu)$. 

Let us recall several notions related to [SSP].

\vskip 2mm\noindent
{\bf [IU]} $\lambda_0$ is an eigenvalue of $L_D$; and there exists,
for any $t>0$,  a constant $C_t >0$ such that
\begin{equation*}
p(x,y,t) \le C_t \  \phi_0(x) \phi_0(y), \quad x,y \in D,
\end{equation*}
where $\phi_0$ is the normalized positive eigenfunction for $\lambda_0$. 

\vskip 2mm\noindent
This notion was introduced by Davies-Simon \cite{DS84}, and investigated
extensively because of its important consequences
(see \cite{Ba91}, \cite{Ba98}, \cite{BD89}, \cite{BD92},  \cite{Dav89}, 
\cite{Me00}, \cite{Me06}, \cite{Mu02}, \cite{Mu07}, \cite{To07},
and references therein). It looks, on the surface, not related  to 
perturbation theory. But it has turned out (\cite{Mu07}) that [IU] implies 
the following condition [SP] for any $a<\lambda_0$.

\vskip 2mm\noindent
{\bf [SP]} The constant function 1 is a small
perturbation of $L-a $ on $D$, i.e.,
for any $\varepsilon > 0$ there exists a compact subset $K$
of $D$ such that 
\begin{equation*}
\int_{D \setminus K} G(x,z) G(z,y) m(z) d\nu(z)
\le \varepsilon G(x,y), \qquad
x,y \in D \setminus K ,
\end{equation*}
where $G$ is the Green function of $L-a $ on $D$.

\vskip 2mm\noindent
This condition is a special case of the notion introduced by 
Pinchover \cite{Pi89}. Recall that [SP] implies the following condition
[SSP] (see \cite{Mu97}).

\vskip 2mm\noindent
{\bf [SSP]} The constant function $1$ is a semismall perturbation of 
$L-a $ on $D$, i.e.,
for any $\varepsilon > 0$ there exists a compact subset $K$
of $D$ such that 
\begin{equation*}
\int_{D \setminus K} G(x^0,z) G(z,y) m(z) d\nu(z)
\le \varepsilon G(x^0,y), \qquad
y \in D \setminus K ,
\end{equation*}
where $x^0$ is a fixed reference point in $D$.

\vskip 2mm\noindent
This condition [SSP] implies that $L_D$ admits a complete orthonormal base 
of eigenfunctions
$\{ \phi_j \}_{j=0}^{\infty}$ with eigenvalues
$\lambda_0 < \lambda_1 \le \lambda_2 \le \cdots$ repeated
according to multiplicity; furthermore, 
for any $j=1,2,\cdots$, the function
$\phi_j /\phi_0$ has a continuous extension $[\phi_j /\phi_0]$
up to the Martin boundary $\partial_M D$ of $D$ for $L-a$
(see Theorem 6.3 of \cite{Pi99}).

\vskip 2mm
We show in this paper that [SSP] also implies the following condition
[SIU].

\vskip 2mm\noindent
{\bf [SIU]} $\lambda_0$ is an eigenvalue of $L_D$; and there exist,
for any $t>0$ and compact subset $K$ of $D$, positive
constants $A$ and $B$ such that
\begin{equation*}
A \  \phi_0(x) \phi_0(y) \le 
p(x,y,t) \le B \  \phi_0(x) \phi_0(y), \quad x \in K,\ \ y \in D.
\end{equation*}

\vskip 2mm\noindent
This notion was introduced by Ba\~nuelos-Davis \cite{BD89}, where they
called it one half  IU. Here we should recall that [IU] implies that 
for any $t>0$ there exists a constant $c_t >0$ such that
\begin{equation*}
c_t \  \phi_0(x) \phi_0(y) \le p(x,y,t), \quad x,y \in D.
\end{equation*}
We see that the same argument as in the proof of Theorem 3.1 in
\cite{Mu93} (or the argument in the proof of Theorem 1.2 below) shows that
[SIU] implies the following condition [NUP] (i.e., non-uniqueness for the positive 
Cauchy problem).

\vskip 2mm\noindent
{\bf [NUP]} The Cauchy problem
\begin{equation} \label{eq:1.3}
(\partial_t + L)u=0 \quad \text{in} \ \  D \times (0,T),
\quad u(x,0) = 0  \quad \text{on} \ \  D
\end{equation}
admits a solution $u$ with $u(x,t)>0$ in $D \times (0,T)$.

\vskip 2mm\noindent
We say that {\bf [UP]} holds for (\ref{eq:1.3}) when 
any nonnegative solution of (\ref{eq:1.3}) is identically zero.
We note that [UP] implies that 
the constant function $1$ is a "big" perturbation of 
$L-a $ on $D$ in some sense (see Theorem 2.1 of  \cite{Mu03}).

Fix $a<\lambda_0$, and suppose that [SSP] holds. Let
$D^* = D \cup \partial_M D$ be the Martin compactification of 
$D$ for $L-a$, which is a compact metric space.
Denote by $\partial_m D$ \  the minimal Martin boundary 
of $D$ for $L-a $, which is a Borel subset of 
the Martin boundary $\partial_M D$
of $D$ for $L-a$.
Here, we note that $\partial_M D$ and $\partial_m D$ are 
independent of $a$ in the following sense: if [SSP] holds,
then for any $b<\lambda_0$ there is a homeomorphism $\Phi$ from 
the Martin compactification of $D$ for $L-a$ onto that for $L-b$
such that $\Phi|_D = identity$, and $\Phi$ maps  
the Martin boundary and minimal Martin boundary of $D$ for $L-a$
onto those for $L-b$, respectively (see Theorem 1.4 of \cite{Mu97}).

Now, we are ready to state our main results. In the following theorems 
we assume that [SSP] holds for some fixed $a<\lambda_0$.

\begin{Thm} \label{th:1.1}
The condition [SSP] implies [SIU].  
\end{Thm}

\begin{Thm} \label{th:1.2}
Assume [SSP]. Then,  
for any $ \xi \in \partial_M D$ there exists the limit
\begin{align} \label{eq:1.4} 
	\lim_{D \ni y \to \xi} 
	\frac{p(x,y,t)}{\phi_0(y)}
	\equiv 
	q(x,\xi,t),\quad
	x\in D,\ t \in {\mathbf{R}}.
\end{align}
Here, as functions of $(x,t)$,
$\{p(x,y,t)/\phi_0(y)\}_y$  converges to $q(x,\xi,t)$ as $ y \to \xi$
uniformly on any compact subset of 
$ D \times {\mathbf{R}}$. 
Furthermore, 
$q(x,\xi,t)$ is a continuous function on 
$D \times \partial_M D \times {\mathbf{R}}$ such that
\begin{equation} \label{eq:1.5}
q>0 \  \  \text{on} \ \  D \times \partial_M D \times(0,\infty ),
\end{equation}
\begin{equation} \label{eq:1.6}
q=0 \  \  \text{on} \ \  D \times \partial_M D \times(-\infty,0],
\end{equation}
\begin{equation} \label{eq:1.7}
(\partial_t + L)q(\cdot,\xi,\cdot)=0 \ \  \text{on} \ \ 
D \times {\mathbf{R}}.
\end{equation}
\end{Thm}

\begin{Thm} \label{th:1.3}
Assume [SSP]. Consider the equation (\ref{eq:1.1}) for 
$I =(0,T)$ with $0 < T \le \infty$.
Then,
for any nonnegative solution $u$ of (\ref{eq:1.1})
there exists a unique pair of Borel measures $\mu$ on $D$ and 
$\lambda$ on $\partial_M D \times [0,T)$
such that
$\lambda$ is supported by the set  $\partial_m D \times [0,T)$, and
\begin{equation} \label{eq:1.8}
u(x,t) = \int_D p(x,y,t) d\mu(y)  +
\int_{\partial_M D \times [0,t)} q(x,\xi,t-s) d\lambda(\xi,s)
\end{equation}
for any $(x,t)\in D \times I$.

Conversely, for any 
Borel measures $\mu$ on $D$ and 
$\lambda$ on $\partial_M D \times [0,T)$
such that
$\lambda$ is supported by $\partial_m D \times [0,T)$
and
\begin{equation} \label{eq:1.9}
\int_D p(x^0,y,t) d\mu(y) < \infty, \quad 0<t<T, 
\end{equation}
\begin{equation} \label{eq:1.10}
\int_{\partial_M D \times [0,t)} q(x^0,\xi,t-s) d\lambda(\xi,s)< \infty ,
\quad 0<t<T, 
\end{equation}
where $x^0$ is a fixed point in $D$, 
the right hand side of (\ref{eq:1.8}) is a
nonnegative solution of (\ref{eq:1.1}) for 
$I =(0,T)$ with $0 < T \le \infty$.
\end{Thm}

The proof of this theorem will be given in Sections 4 and 5.
It is based upon the abstract integral representation theorem 
established in \cite{Mu07}, without assuming [IU], via a
parabolic Martin representation theorem and Choquet's theorem 
(see \cite{Jan71}, \cite{Ma91},\cite{Ph66}).
Its key step is to identify the parabolic Martin boundary.

This theorem is an improvement of Theorem 1.2 of \cite{Mu07};
where the condition [IU], which is more stringent than [SSP],
is assumed. It is also an answer to a problem raised in Remark 4.13 
of \cite{Mu07}. Note that (\ref{eq:1.8}) gives explicit integral representations
of nonnegative solutions to (\ref{eq:1.1})  provided that the Martin boundary
$\partial_M D$ of $D$ for $L-a$ is determined explicitly.
We consider that [SSP] is one of the weakest possible condition for 
getting such explicit integral representations.

Let us recall that when [UP] hods for (\ref{eq:1.3}), the structure of all
nonnegative solutions to (\ref{eq:1.1}) for $I =(0,T)$ is extremely
simple. Namely, the following theorem holds (see \cite{AT92}).

\vskip 2mm\noindent
{\bf Fact AT}  Assume [UP]. Then,
for any nonnegative solution $u$ of (\ref{eq:1.1})
with $I =(0,T)$,
there exists a unique Borel measure $\mu$ on $D$ 
such that
\begin{equation} \label{eq:1.11}
u(x,t) = \int_D p(x,y,t) d\mu(y),
\quad  (x,t) \in D \times I.
\end{equation}
Conversely, for any 
Borel measure $\mu$ on $D$ satisfying (\ref{eq:1.9}),
the right hand side of (\ref{eq:1.11}) is a nonnegative
solution of (\ref{eq:1.1}) with $I =(0,T)$.

\vskip 2mm
It is quite interesting that when [UP] holds, the elliptic Martin 
boundary disappears in the parabolic representation theorem; while it 
enters in many cases of [NUP].

Finally, we state an integral representation 
theorem for the case  $I = (-\infty,0)$.

\begin{Thm} \label{th:1.4}
Assume [SSP]. 
Consider the equation (\ref{eq:1.1}) for $I = (-\infty,0)$.
Then,
for any nonnegative solution $u$ of (\ref{eq:1.1})
there exists a unique pair of 
a nonnegative constant $\alpha$ and a
Borel measure $\lambda$ on $\partial_M D \times (-\infty,0)$
supported by the set  $\partial_m D \times (-\infty,0)$
such that
\begin{equation} \label{eq:1.12}
u(x,t) = \alpha e^{-\lambda_0 t} \phi_0(x) +
\int_{\partial_M D \times (-\infty,t)} q(x,\xi,t-s) d\lambda(\xi,s) 
\end{equation}
for any $(x,t)\in D \times (-\infty,0)$. 

Conversely, for any 
nonnegative constant $\alpha$ and a
Borel measure $\lambda$ on $\partial_M D \times (-\infty,0)$
such that it is supported by $\partial_m D \times (-\infty,0)$
and
\begin{equation} \label{eq:1.13}
\int_{\partial_M D \times (-\infty,t)} q(x^0,\xi,t-s) d\lambda(\xi,s)< \infty ,
\quad -\infty <t<0, 
\end{equation}
the right hand side of (\ref{eq:1.12}) is a
nonnegative solution of (\ref{eq:1.1}).
\end{Thm}

This theorem is an improvement of Theorem 6.1 of 
\cite{Mu07}, where [IU] is assumed instead of [SSP].

Here, in order to illustrate a scope of Theorems \ref{th:1.3}
and \ref{th:1.4}, we give a simple example.
Further examples will be given in Section 7.

\begin{Ex} \label{ex:1.5}
Let $D$ be a domain in $\mathbf{R}^2$ with
finite area. Then, by Theorem 6.1 of \cite{Mu05}, the
constant function $1$ is a small perturbation of 
$L = - \Delta$ on $D$. Thus Theorems \ref{th:1.3}
and \ref{th:1.4} hold true for the heat equation 
\begin{equation*} 
(\partial_t -\Delta)u=0 \quad \text{in} \quad D \times I. 
\end{equation*}
Note that there exist many bounded planar domains for which the heat
semigroup is not intrinsically ultracontractive (see Example 1 of 
\cite{DS84} and Section 4 of \cite{BD89}). Thus, the last assertion 
of this example is new for such domains.
\end{Ex}

The remainder of this paper is organized as follows.
In Section 2 we prove Theorem~\ref{th:1.1}, and Theorem~\ref{th:1.2}
is proved in Section 3. Sections 4 and 5 are devoted to the proof of
Theorem~\ref{th:1.3}. In Section 4 we show it in the case of 
$I = (0,\infty)$. In Section 5 we show it in the case of 
$I = (0,T)$ with $0 < T <\infty$ by making use of results to be given
in Section 4. Theorem~\ref{th:1.4} is proved in Section 6.
Finally we shall give two more concrete examples in Section 7 with
emphasis on sharpness of concrete sufficient conditions of [SSP].

\section{[SSP] implies [SIU]} \label{sec:2}

In this section we prove Theorem~\ref{th:1.1}.

\vskip 2mm
\noindent{\bf Proof of Theorem~\ref{th:1.1}} \quad
We may and shall assume that
$ a = 0 <\lambda_0$.
Let $ G$ be the Green function of $ L$ on $D$.
For any $t >0$, put 
\begin{align*} 
	G_t(x,y) & = 
	\int_t^\infty p(x,y,s)\,ds,\\
	G^t(x,y) & = 
	\int_0^t p(x,y,s)\,ds.
\end{align*}
Then $G= G_t+G^t$.
Let us show that for any $t>0$ and any compact subset $K$ of $D$
there exists a constant $A>0$ such that
\begin{equation} \label{eq:2.1}
A \,  \phi_0(x) \,\phi_0(y)\, \leq \,p(x,y,t), \quad x \in K,\ \ y \in D.
\end{equation}
Fix a compact subset $K$. We may assume that $x^0 \in K$.  
Let $K_1 \subset D$ be a compact neighborhood of $K$. Then 
the same argument as in the proof of Theorem 1.5 of \cite{Mu97}
shows that 
\begin{equation} \label{eq:2.2} 
C^{-1}\, G(x^0,z) \le \phi_0(z) \le C \,G(x^0,z), \quad 
z \in D \setminus K_1,
\end{equation}
for some constant $C>0$. Fix $t>0$, and put
\[
 \epsilon_{t}= \frac{1}{2 \lambda_0}\, \left(\,1-e^{-t\lambda_0}\,\right).
\]
By [SSP] and (\ref{eq:2.2}), there exits a compact subset 
 $K_2 \supset K_1$  such that 
\begin{equation}\label{eq:2.3}
 \int_{D \setminus K_2} \phi_0 (z)\,G(z,y)\,d\mu(z)\, 
\le \epsilon_{t}\, \phi_0 (y),
 \quad y \in D \setminus K_2,
\end{equation}
where $d\mu(z)=m(z)\,d\nu(z)$. Since 
\[ \frac{\phi_0 (y)}{\lambda_0} = \int_{D} G(y,z)\,\phi_0 (z)\,d\mu(z),\]
and  $G(y,z)=G(z,y)$, (\ref{eq:2.3}) yields 
\begin{eqnarray}\nonumber
\frac{\phi_0 (y)}{\lambda_0} \leq  \int_{K_2} G_t(z,y)\,\phi_0 (z)\,d\mu(z)
&+& \int_{K_2} G^t(z,y)\,\phi_0 (z)\,d\mu(z) \nonumber\\
&+&  \epsilon_{t} \, \phi_0(y)  \label{eq:2.4}
\end{eqnarray}
for any $ y \in D \setminus K_2$. By Fubini's theorem,  
\begin{eqnarray*}
\int_{D} G_{t}(z,y)\,\phi_0 (z)\,d\mu(z)
&=&  \,\int_{t}^{\infty}\,ds\,\int_{D} p(z,y,s)\,\phi_0 (z)\,d\mu(z) \\
&=&\,\int_{t}^{\infty}\,e^{-\lambda_0 s}\,\phi_0 (y)\,ds\\
&=&\frac{1}{\lambda_0}\, e^{-\lambda_0 t}\,\phi_0 (y).
\end{eqnarray*}
Thus 
\[ 
\int_{K_2} G_{t}(z,y)\,\phi_0 (z)\,d\mu(z)\,\leq\,\frac{1}{\lambda_0}\, 
e^{-\lambda_0 t}\,\phi_0 (y).
\]
This together with (\ref{eq:2.4}) implies 
\begin{equation}\label{eq:2.5}
\epsilon_t \,\phi_0 (y)\,\leq\, \int_{K_2}  G^t(z,y)\,\phi_0 (z)\,d\mu(z) .
\end{equation}
Choose a compact subset $K_3$ whose interior includes $K_2$. By the 
parabolic Harnack inequality, there exists a constant $C_1$ depending on
$t, K_2, K_3$ such that 
\[
p(z,y,s)\leq C_1 \,p(x,y,2t),
\] 
for any $ x,z \in K_2$, $y \in D\setminus K_3$, and $0<s\leq t$.
We have 
\begin{eqnarray}\nonumber
G^t(z,y) &=& \int_{0}^{t}\, p(z,y,s)\,ds \\
&\leq& C_1 \, t \, p(x^0,y,2t),  \qquad
z \in K_2, \  y \in D\setminus K_3. \label{eq:2.6}
\end{eqnarray}
Thus
\[
\int_{K_2} G^{t}(z,y)\,\phi_0 (z)\,d\mu(z) \leq \,\left[\,C_1\,t\,\int_{K_2} 
\phi_0 (z)\,dz\right]
\, p(x^0,y,2t).
\]
This together with (\ref{eq:2.5}) implies 
\begin{equation}\label{eq:2.7}
\phi_0 (y)\, \leq \,C_2\, p(x^0,y,2t), \quad y \in D\setminus K_3,
\end{equation}
where $$C_2=\frac{1}{\epsilon_t}\,C_1\,t \int_{K_2} \phi_0(z)d\mu(z).$$ 
By the parabolic Harnack inequality, 
\[
p(x^0,y,2t) \leq C \, p(x,y,3t),
\quad x \in K, \, y \in D,
\]
for some constant $C>0$. 
This together with (\ref{eq:2.7}) yields 
the desired inequality (\ref{eq:2.1}). 
It remains to show that for any $t>0$ and a compact subset $K$ of $D$  
there exists a constant $B$ such that 
\begin{equation}\label{eq:2.8}
p(x,y,t)\,\leq\, B \, \phi_0(x)\,\phi_0 (y),\quad  x \in K, \, y \in D.
\end{equation}
Fix a compact subset $K$. We may assume that $x^0 \in K$. 
Let $K_1 \subset D$ be a compact neighborhood of $K$. By the parabolic 
Harnack inequality
there exists a constant $c>0$ such that
\[
c \,p(x^0,y,t) \leq p(z,y,2t), \quad z \in K_1, \  y \in D.
\]
Thus, for any $y \in D$, 
\begin{eqnarray}\nonumber
e^{-2t\lambda_0\,}\,\phi_0 (y) &=& \int_D \,\phi_0(z)\,p(z,y,2t)\,d\mu(z)\\\nonumber
&\geq& \int_{K_1} \,\phi_0(z)\,p(z,y,2t)\,d\mu(z)\\\nonumber
&\geq& \,c\,\left[\,\int_{K_1} \,\phi_0(z)d\mu(z)\,\right]\,p(x^0,y,t).\,
\end{eqnarray}
This implies (\ref{eq:2.8}), since
\[
C \, p(x^0,y,t) \geq p(x,y,t/2),
\quad x \in K, \, y \in D,
\]
for some constant $C>0$. 
(We should note that in proving  (\ref{eq:2.8}) 
we have only used the consequence of [SSP] that $\phi_0$ is a 
positive eigenfunction.)
$\quad\square$

\begin{Rem} \label{rem:2.1}
It is an open problem whether [SIU] implies [SSP] or not.
Furthermore, the problem whether [SSP] implies [SP] or not in the case $n>1$ is 
still open.
\end{Rem} 

                                                                                                                       
\section{Parabolic Martin kernels} \label{sec:3}

In this section we prove Theorem~\ref{th:1.2}. Throughout the present 
section we assume [SSP]. We may and shall assume that $a=0 < \lambda_0$.  
Let $G$ be the Green function of $L$ on $D$. For any $0< \delta < t$, put 
\begin{equation} \label{eq:3.1}
G_{\delta}^t(x,y)=\int_{\delta}^t\,p(x,y,s)\,ds. 
\end{equation}
We denote by $\partial_{M} D$ the Martin boundary of $D$ for $L$. In order to 
prove Theorem~\ref{th:1.2}, we need two lemmas. 

\begin{Lem}\label{lm:3.1}
Let  $\xi \in \partial_{M} D$. Suppose that a sequence 
$\{y_n\}_{n=1}^\infty \subset D$  converges to   $\xi$, and  
there exists the limit
\begin{equation}\label{eq:3.2}
\lim_{n \to \infty} \frac{G_{\delta}^t(z,y_n)}{\phi_0(y_n)}=w(z,t),\quad z \in D.
\end{equation}
Then  
\begin{equation}\label{eq:3.3}
\lim_{n \to \infty}  \int_{D} \, G(x,z)\,
\frac{G_{\delta}^t(z,y_n)}{\phi_0(y_n)}\,\,d\mu(z)=  
\int_{D} \, \,G(x,z)\,w(z,t) \,d\mu(z)
\end{equation}
for any  $x \in D$, where $d\mu(z)=m(z)d\nu(z)$.
\end{Lem}

\vskip 2mm
\noindent{\bf Proof}\quad
Fix $x \in D$. Let   $K_1 \subset D$ be a 
compact neighborhood of $x$. By [SSP], there exists a constant 
$C>0$ such that  
\begin{equation} \label{eq:3.4}
C^{-1}\, \phi_0 (y)\, \leq G(x,y)\, \leq \,C\,\phi_0(y),\quad 
y \in  D \setminus K_1.
\end{equation}
Let $\epsilon >0$. Then there exists a compact subset  $K \supset K_1 $ such that  
\begin{equation} \nonumber
\int_{D \setminus K} G(x,z)\,\frac{G(z,y)}{G(x,y)}\,d\mu(z) < 
\frac{\epsilon}{3 C}, \quad y \in D\setminus K.
\end{equation}
Thus, for $n$ sufficiently large, 
\begin{eqnarray*}
 \int_{D \setminus K} G(x,z)\,\left[\,
\frac{G_{\delta}^t(z,y_n)}{\phi_0(y_n)}\,\right]\,d\mu(z)\, &\leq& \,
 \int_{D \setminus K} G(x,z)\,\left[\,\frac{C\,G(z,y_n)}{G(x,y_n)}\,\right]\,
d\mu(z)\\
& <& \frac{\epsilon}{3}.
\end{eqnarray*}
By Fatou's lemma,  
\begin{equation}\nonumber
\int_{D \setminus K} G(x,z)\,w(z,t)\,d\mu(z)\, \leq \,\frac{\epsilon}{3}.
\end{equation}
By Theorem~\ref{th:1.1}, there exist
constants $A_1$ and $A_2$ such that
\begin{equation*}
A_1 \  \phi_0(x) \phi_0(y) \le 
p(x,y,\delta) \le A_2 \  \phi_0(x) \phi_0(y), \quad x \in K,\ \ y \in D.
\end{equation*}
Then, for any $t > \delta$, the semigroup property yields
\begin{equation}\label{eq:3.5}
A_1\, e^{-\lambda_0 (t-\delta)} \,  \phi_0(x) \phi_0(y) \le 
p(x,y,t) \le \,A_2\,e^{-\lambda_0 (t-\delta)}\,  \phi_0(x) \phi_0(y) 
\end{equation}
for any $x \in K, \,  y \in D$. 
Thus there exists a constant  $B>0$ such that  for any $n$
\begin{equation*}
\frac{G_{\delta}^{t}(z,y_n)}{ \phi_0 (y_n)} \leq  B\,\phi_0 (z),\quad z \in K.
\end{equation*}
Then Lebesgue's dominated convergence theorem yields
\begin{equation}\nonumber
\lim_{n \to \infty}  \int_{ K} \,G(x,z)\,\left[\, 
\frac{G_{\delta}^{t}(z,y_n)}{\phi_0(y_n)}\,\right]\,d\mu(z)\,
= \,\int_{K}\,G(x,z)\,w(z,t) \,d\mu(z).
\end{equation}
Therefore, for $n$ sufficiently large, 
\[
\left|\, \int_{ D} \,G(x,z)\,\left[\, 
\frac{G_{\delta}^{t}(z,y_n)}{\phi_0(y_n)}\,\right]\,d\mu(z)\,- \,
\int_{D}\,G(x,z)\,w(z,t) \,d\mu(z)\,\right| < \epsilon.
\]
This shows (\ref{eq:3.3}).$\quad\square$

By Lemma 6.1 of \cite{Pi99}, it follows from [SSP] that there exists 
the limit 
\begin{equation}\label{eq:3.6}
\lim_{D \ni y \to \xi} \frac{G_{D}(y,z)}{\phi_0(y)}=h(\xi,z),\quad 
(\xi,z) \in \partial_{M}D \times D,
\end{equation}
and $h$ is a positive continuous function on $\partial_M D \times D$. 
From this we show the following lemma. 
\begin{Lem}\label{lm:3.2}
Under the same assumptions as in Lemma~\ref{lm:3.1}, one has 
\begin{eqnarray}\nonumber
\int_{D} h(\xi,z)\,G_{\delta}^t(z,x)\,d\mu(z)&=& \lim_{n \to \infty}  
\int_{D} \, \frac{G(y_n,z)}{\phi_0(y_n)}\,G_{\delta}^t(z,x)\,d\mu(z)\\
&=& \int_{D} \, G(x,z)\,w(z,t)\,d\mu(z)  \label{eq:3.7}
\end{eqnarray}
for any $ x \in D$.
\end{Lem}

\vskip 2mm
\noindent{\bf Proof } Fix $x \in D$. Let $K_1 \subset D$  be a compact 
neighborhood of $x$. 
By Theorem~\ref{th:1.1},  (\ref{eq:3.4}) and (\ref{eq:3.5}), there exists a 
constant $C_1 >0$ such that 
\[
C_1\,G(z,x)\,\leq\,G_{\delta}^t(z,x)\,\leq\,\,G(z,x),\quad 
z \in D\setminus K_1.
\]
Let $\epsilon > 0$. By [SSP], there exists a compact subset $K \supset K_1$ 
such that  
\begin{equation}\label{eq:3.8}
 \int_{D \setminus K} \,
\left[\,\frac{G(y_n,z)}{\phi_0(y_n)}\,\right]\,G_{\delta}^t(z,x)\,d\mu(z)\, 
< \frac{\epsilon}{3},
\end{equation}
for $n$ sufficiently large. By Fatou's lemma,
\begin{equation}\label{eq:3.9}
 \int_{D \setminus K} \,h(\xi,z)\,\,G_{\delta}^t(z,x)\,d\mu(z)\, 
\leq\, \frac{\epsilon}{3}.
\end{equation}
On the other hand, for any sufficiently large $n$
\[
 \left[\,\frac{G(y_n,z)}{\phi_0(y_n)}\,\right]\,G_{\delta}^t(z,x)\, 
\leq\, C_2,\quad z \in K ,
\]
where $C_2$ is a positive constant. By Lebesgue's dominated convergence theorem,
\begin{equation}\label{eq:3.10}
\lim_{n \to \infty}  \int_{ K} \, 
\frac{G(y_n,z)}{\phi_0(y_n)}\,G_{\delta}^t(z,x)\,d\mu(z)\,
= \,\int_{K} h(\xi,z)\,G_{\delta}^t(z,x)\,d\mu(z).
\end{equation}
Combining (\ref{eq:3.8}), (\ref{eq:3.9}) and (\ref{eq:3.10}), we get 
the first equality. It remains to show the second equality of (\ref{eq:3.7}).
By Fubini's theorem and the symmetry 
\[p(x,y,t)=p(y,x,t),\]
we have 
\begin{eqnarray*} 
\int_{D} \, G(y_n,z)\,G_{\delta}^t(z,x)\,d\mu(z)&=& \int_{0}^{\infty} \,dr\,
\int_{\delta}^{t} \,ds\,p(y_n,x,r+s)\\
&=& \int_{D} \, G(x,z)\,G_{\delta}^t(z,y_n)\,d\mu(z).
\end{eqnarray*}
This together with Lemma~\ref{lm:3.1} implies the second equality.$\quad\square$

\vskip 2mm
\noindent{\bf Proof of Theorem 1.2 } 
Let $\{y_j\}_{j=1}^\infty \subset D$ be any  sequence converging to 
$\xi \in \partial_M D$. Put 
\begin{equation} \label{eq:3.11}
u_{j}(x,t)=\frac{p(x,y_j,t)}{\phi_0(y_j)} \quad \text{ for } t>0, \qquad
u_{j}(x,t)=0 \quad \text{ for } t \leq 0.
\end{equation}
Since [SIU] holds, it follows from 
the parabolic Harnack inequality and local a priori estimates for 
nonnegative solutions to parabolic equations 
(see \cite{Aro68} and \cite{Is99}) that there exists a 
subsequence $\{u_{j_k}\}_{k=1}^\infty$ such that $u_{j_k}$ converges, 
as $k \to \infty$, uniformly 
on any compact subset of $D\times \mathbf{R}$ to a solution $u$ of the equation 
\[
\left(\,\partial_t+L\,\right)\,u\,=\,0 \  \text{ in } D\times \mathbf{R}
\]
satisfying $ u>0$ on $ D\times (0,\infty)$ and $u=0$ on $ D\times (-\infty,0]$. 
Thus, in order to prove Theorem~\ref{th:1.2}, it suffices to show that 
the limit function $u$ is independent of $\{y_{j_k}\}_{k=1}^\infty$ and 
uniquely determined by $\xi$.
Let  $\{y_j\}_{n=1}^\infty$ and $\{y_j'\}_{n=1}^\infty$ be two sequences in $D$ 
converging to $\xi$. Define $u_j$ by (\ref{eq:3.11}), and $u_j'$ by 
(\ref{eq:3.11}) with $y_j$ replaced by $y_j'$. Suppose that $\{u_j\}_{j=1}^\infty$ 
and  $\{u_j'\}_{j=1}^\infty$ converge to $u$ and $u'$, respectively.
For any $t>\delta>0$, put 
\[
w(z,t)=\int_{\delta}^t\,u(z,s)\,ds, \qquad  
w'(z,t)=\int_{\delta}^t\,u'(z,s)\,ds.
\]
Then we have 
\[
 \lim_{n \to \infty} \frac{G_{\delta}^{t}(z,y_n)}{\phi_0(y_n)}=w(z,t), \qquad
\lim_{n \to \infty} \frac{G_{\delta}^{t}(z,y_n')}{\phi_0(y_n')}=w'(z,t).
\]
By Lemma~\ref{lm:3.2},
\begin{eqnarray*}
\int_{D} \, G(x,z)\, w(z,t)\,d\mu(z)&=&\int_{D} \, h(\xi,z)\, 
G_{\delta}^t(z,x)\,d\mu(z)\\
&=&\int_{D}\,G(x,z)\,  w'(z,t)\,d\mu(z).
\end{eqnarray*}
Thus $w(x,t)=w'(x,t)$, which implies $u(x,t)=u'(x,t)$. This completes 
the proof of Theorem~\ref{th:1.2}.
\qed

\section{Integral representations; the case $I=(0,\infty)$} \label{sec:4}

In this section we prove Theorem~\ref{th:1.3} in the case  $T=\infty$.

We first state an abstract integral representation theorem which holds 
without [SSP]. 
For $x \in D$ and $ r>0$, we denote by $B(x,r)$ the geodesic ball in the
Riemannian manifold $M$ with center $x$ and radius $r$.
Let $ x^0$ be a reference point in $D$.
Choose a nonnegative continuous function $a$ on $D$ such that
$a(x)=1$ on $B(x^0,r^0)$ and $a(x)=0$ outside $B(x^0,2r^0)$ 
for some $ r^0 >0$ with $B(x^0,3r^0) \Subset D$.
Choose a nonnegative continuous function $b$ on $\mathbf{R}$
such that 
$0< b(t)< e^{\gamma t}$ on $(1,\infty)$ for some $\gamma < \lambda_0$,
and $b(t)=0$ on $(-\infty,1]$.
Denote by $\beta$ the measure defined by 
$ d \beta(x,t)= a(x)b(t)m(x)\,d \nu(x) dt $.
For any nonnegative measurable function $u$ on $ Q = D \times (0,\infty)$,
we write 
\[
	\beta(u)= \iint_{Q}\, u(x,t)\, d\beta(x,t).
\]
Denote by $ P(Q)$ the set of all nonnegative solutions of \rfe{eq:1.1}
with $I=(0,\infty)$, and put 
\[
	P_{\beta}(Q)=\left\{
	u \in P(Q) ; \beta(u)< \infty
	\right\}.
\]
Note that for any $ u \in P(Q)$ there exists a function $ b$ as above 
such that $ \beta(u) < \infty$; thus 
$ P(Q)= \bigcup_{\beta} P_{\beta}(Q)$.
Furthermore, the parabolic Harnack inequality shows that if $ \beta(u)=0$,
then $u=0$.
Now, let us define the $\beta$-Martin boundary $ \partial^\beta_M Q$ of $Q$
with respect to $ \partial_t +L$
along the line given in \cite{Ma91} and \cite{Jan71}. Put
\begin{align*}
	\begin{array}{ll}
	 p(x,t ; y,s)= p(x,y,t-s), & t>s,\ x,y \in D,\\
	 \noalign{\vskip 2mm}
	 p(x,t ; y,s)=0,           & t \leq s ,\ x,y \in D.
	 \end{array}
\end{align*}
Define the $\beta$-Martin kernel $K_\beta$ by 
\[
	K_\beta(x,t;y,s)=
	\frac{p(x,t;y,s)}{\beta\left( p(\,\cdot\,;y,s)\right)},
	\quad (x,t),\ (y,s) \in Q,
\]
where $ \beta\left( p(\,\cdot\,;y,s)\right)= 
\iint_{Q}\, p(z,r ; y,s)\,
d \beta(z,r)  $.
Note that $\beta\left( p(\,\cdot\,;y,s)\right) < \infty$
for any $(y,s) \in Q$, since 
$0< b(t)< e^{\gamma t}$ on $(1,\infty)$ for some $\gamma < \lambda_0$.
Let $ \{ D_j\}_{j=1}^\infty$ be an exhaustion of $D$ such that
each $D_j$ is a domain with smooth boundary, $ D_j \Subset D_{j+1} 
\Subset D$, $ \bigcup_{j=1}^\infty D_j = D$, and $ B(x^0,3r^0) \Subset D_1$.
Put $Q_j=D_j \times (1/j,j)$.
For $ Y=(y,s),\ Z=(z,r) \in Q$, let 
\[
	\delta_\beta(Y,Z)
	= \sum_{j=1}^\infty 2^{-j} \sup_{X \in Q_j}
	\frac{\left| K_\beta(X;Y) - K_\beta(X;Z)\right|}
	{1+\left| K_\beta(X;Y) - K_\beta(X;Z)\right|}.
\] 
Then we see that $ \delta_\beta$ is a metric on $Q$, and the topology 
on $Q$ induced by $\delta_\beta $ is equivalent to the original 
topology of $Q$. Denote by $ Q^{\beta*}$ the completion of $Q$ with respect to 
the metric $\delta_\beta$.
Put $ \partial^\beta_M Q = Q^{\beta*} \setminus Q$.
A sequence $ \{ Y^k\}_{k=1}^\infty$ in $Q$ is called a fundamental sequence if 
$ \{ Y^k\}_{k=1}^\infty$ has no point of accumulation in $Q$ and 
$ \left\{K_\beta(\,\cdot \, ;  Y^k ) \right\}_{k=1}^\infty$
converges uniformly on any compact subset of $Q$ to a nonnegative 
solution of \rfe{eq:1.1} with $I=(0,\infty)$.
By the local a priori estimates for solutions of \rfe{eq:1.1},
for any $ \Xi \in \partial^\beta_M Q $ there exist a unique nonnegative solution
$ K_\beta(\,\cdot\, ; \Xi)$ of \rfe{eq:1.1} and a fundamental sequence 
$ \{ Y^k\}_{k=1}^\infty$ in $Q$ such that
\[
	\lim_{k \to \infty}\sum_{j=1}^\infty 2^{-j} \sup_{X \in Q_j}
	\frac{\left| K_\beta(X;Y^k) - K_\beta(X;\Xi)\right|}
	{1+\left| K_\beta(X;Y^k) - K_\beta(X;\Xi)\right|}=0.
\] 
Thus the metric $\delta_\beta$ is canonically extended to $Q^{\beta*}$.
Furthermore, $ Q^{\beta*}$ becomes a compact metric space,
since by the parabolic Harnack inequality,
any sequence $ \{ Y^k\}_{k=1}^\infty$ with no point of 
accumulation in $Q$ has a fundamental subsequence.  We call 
$ K_\beta(\,\cdot\, ;\Xi)$, $ \partial^\beta_M Q $ and $ Q^{\beta*}$ 
the $\beta$-Martin kernel, $\beta$-Martin boundary and 
$\beta$-Martin compactification for $(Q,\ \partial_t + L)$, respectively.
Note that $\beta \left( K_\beta( \,\cdot\, ; \Xi) \right)\leq 1$
by  Fatou's lemma; and so $ K_\beta(\,\cdot\, ; \Xi) \in P_{\beta}(Q)$.
A nonnegative solution $ u \in P_{\beta}(Q)$ is said to be minimal
if for any nonnegative solution $ v \leq u$ there exists a nonnegative 
constant $ C$ such that $ v=Cu$. Put
\[
	\partial^\beta_m Q= \left\{
	\Xi \in \partial^\beta_M Q ; K_\beta(\,\cdot\, ; \Xi)\ 
	\mbox{is minimal and}\ \beta \left( K_\beta( \,\cdot\, ; \Xi) \right)
	= 1
	\right\},
\]
which we call the minimal $ \beta$-Martin boundary for  $(Q,\ \partial_t + L)$.

Observe that $ D \times [0,\infty)$ is embedded into $Q^{\beta*}$, and 
$D \times \{ 0\} \subset \partial^\beta_M Q$.
Indeed, with $ y \in D$ fixed, for any sequence $ \{ Y^k\}_{k=1}^\infty$ in $Q$
with $\lim_{k \to \infty} Y^k=(y,0)$ we have 
$\lim_{k \to \infty} K_\beta(x,t ; Y^k)=
p(x,t ; y,0)/\beta\left(p(\,\cdot\, ; y,0) \right)$ ; 
furthermore,  
$ K_\beta(\,\cdot\,  ; y,0) \ne K_\beta(\,\cdot\,  ; z,0)$ if $ y \ne z$.
We also note that any sequence $ \left\{ Y^k =\right.$ 
$ \left. (y^k, s^k )\right\}_{k=1}^\infty$
in $Q$ with $ \lim_{k \to \infty} s^k =\infty$ is a fundamental sequence, 
since $ \lim_{k \to \infty} K_\beta(\,\cdot\,  ; Y^k)=0$.
We denote by $\varpi $ the point in $\partial^\beta_M Q$
corresponding to the Martin kernel which is identically zero : 
$ K_\beta(\,\cdot\, ; \varpi)=0$.
Put 
\[
	{\cal L}^\beta_m Q = \partial^\beta_m Q \setminus
	\left( D \times \{0\} \cup \{ \varpi \}\right).
\]

We obtain the following abstract integral representation theorem 
in the same way as in the proof of Theorem 2.1 and Lemma 2.2 of 
\cite{Mu07}.

\begin{Thm} \label{th:4.1}
For any $ u \in P_\beta(Q)$, there exists a unique pair of finite
Borel measures $ \kappa$ on $D$ and $\lambda$ on 
$ \partial^\beta_M Q \setminus ( D \times \{ 0\})$ such that
$ \lambda$ is supported by the set $ {\cal L}^\beta_m Q$,
\begin{align} \label{eq:4.1}  
	u(x,t)= \int_D \frac{p(x,t ; y,0)}{\beta\left(p( \,\cdot\, ; y,0) 
	\right)}\, d \kappa(y)+
	\int_{{\cal L}^\beta_m Q}K_\beta(x,t ; \Xi)\,d \lambda(\Xi)
\end{align}
for any $(x,t)\in Q$, and 
\begin{align} \label{eq:4.2}  
	\beta(u)= \kappa(D) + \lambda( {\cal L}^\beta_m Q).
\end{align}
Furthermore, the function
\begin{equation*}
v(x,t)= u(x,t) -  \int_D \frac{p(x,t ; y,0)}{\beta\left(p( \,\cdot\, ; y,0) 
\right)}\, d \kappa(y)
\end{equation*}
is a nonnegative solution of the equation
\begin{equation*} 
(\partial_t + L)v=0 \ \  \text{in} \ \ 
D \times {\mathbf{R}}
\end{equation*}
such that $v=0$ on $D \times (-\infty,0]$.

Conversely, for any finite Borel measures 
$ \kappa$ on $D$ and $\lambda$ on 
$ \partial^\beta_M Q \setminus ( D \times \{ 0\})$ such that
$ \lambda$ is supported by the set $ {\cal L}^\beta_m Q$,
the right hand side of \rfe{eq:4.1} belongs to $ P_\beta(Q)$.
\end{Thm}

We put
\[
	P^0_\beta(Q)=
	\left\{v \in P_\beta(Q) ; \lim_{t \downarrow 0} v(x,t)=0\ 
	\mbox{on}\ D
	\right\}.
\]
We show Theorem~\ref{th:1.3} on the basis of Theorem~\ref{th:4.1}.
To this end it suffices to show (\ref{eq:1.8}) for 
$u \in P^0_\beta(Q)$. The key step in the proof is to identify 
${\cal L}^\beta_m Q$. Under the condition [SSP], we shall show that
${\cal L}^\beta_m Q = \partial_m D \times [0,\infty)$. In the
remainder of this section we assume [SSP]. 
We may and shall assume that $ a = 0 <\lambda_0$.

\begin{Lem}\label{lem:4.2}
For any domains $U$ and $W$ with $ U \Subset W \Subset D$,
there exist positive constants $C$ and $\alpha$ such that
\begin{align} \label{eq:4.3} 
	& p(x,y,t) \leq C f(t) \phi_0(x)\phi_0(y),
	\quad x \in U,\ y \in D \setminus W,\ t>0, 
\end{align}
where $ f(t)=e^{-\alpha/t} $ for $0<t<1$, and 
$ f(t)= e^{- \lambda_0 t}$ for $ t \geq 1$.
Furthermore,
\begin{align} \label{eq:4.4} 
	 q(x,\xi,t) & \leq C  f(t) \phi_0(x),
	\quad x \in U,\ \xi \in \partial_M D, \ t>0, \\
	 G(x,y) & \leq C \phi_0(x)\phi_0(y),
	\quad x \in U,\ y \in D \setminus W,  \label{eq:4.5} 
\end{align}
where $G$ is the Green function of $L$ on $D$.
\end{Lem}

This lemma is shown in the same way as Lemmas 4.2 and 4.4 of 
\cite{Mu07}.

Let $K(x,\xi)$ be the Martin kernel for  $L$ on $D$ with reference
point $ x^0 \in D$, i.e., $K(x^0,\xi)=1,\ \xi \in \partial_M D$.
The following lemma gives a relation between $K$ and $q$.

\begin{Lem}\label{lem:4.3}
For any $ \xi \in \partial_M D$,
\begin{align} \label{eq:4.6} 
	 \lim_{D \ni y \to \xi}
	\frac{G(x,y)}{\phi_0(y)}
	&= \int_0^\infty q(x,\xi,t)\,dt,\quad
	x \in D,\\
	 K(x,\xi)&=
	\frac{\int_0^\infty q(x,\xi,t)\,dt}
	{\int_0^\infty q(x^0,\xi,t)\,dt},\quad
	x \in D. \label{eq:4.7}
\end{align}
\end{Lem}

This lemma is shown in the same way as Lemma 4.5 of \cite{Mu07}

\begin{Lem}\label{lem:4.4}
Let $\xi,\eta \in \partial_M D$, $ 0 \leq s,r < \infty$ and $ C >0$. 
If
\[
	q(x,\xi,t-s)=Cq(x,\eta,t-r),\quad (x,t)\in Q,
\]
then $\xi=\eta,\ s=r$ and $C=1$.
\end{Lem}

\noindent{\bf Proof}\quad
Since $q(x,\xi,\tau)>0$ for $\tau>0$ and $q(x,\xi,\tau)=0$
for $\tau \le 0$, we obtain that $s=r$. Thus
$q(x,\xi,\tau)=q(x,\eta,\tau)$. This together with \rfe{eq:4.7}
implies that $K(\,\cdot\, ,\xi)=K(\,\cdot\, ,\eta)$
on $D$. Hence $\xi=\eta$, and so $C=1$.  \qed

Now, let $\beta$ be a measure on $ Q = D \times (0,\infty)$ as 
described in the beginning of this section:
$d  \beta(x,t)= a(x) b(t) m(x) \,d\nu(x)\,dt$.
The following proposition determines the $\beta$-Martin boundary 
$\partial_M^\beta Q$, $\beta$-Martin compactification $ Q^{\beta*}$,
and $\beta$-Martin kernel $K_\beta$ for $ \left( \partial_t+L, Q \right)$.
Recall that
$ p(x,t ; y,s)=p(x,y,t-s)$ and $ K_\beta(\,\cdot\, ;y,s)
= p(\,\cdot\,  ; y,s)/\beta \left( p(\,\cdot\,  ; y,s) \right)$.
We write
\[
	q(x,t;\xi,s)=q(x,\xi,t-s)
\]
for $ \xi \in \partial_M D$ and $ 0 \leq s <\infty$.

\begin{Prop}\label{pro:4.5}
(i)\ The $\beta$-Martin boundary $\partial_M^\beta Q$ of $Q$ for 
$ \partial_t +L$ is equal to the disjoint union of $D \times \{0\}$,
$ \partial_M D\times [0,\infty)$ and the one point set $ \{ \varpi \}$:
\begin{align} \label{eq:4.8} 
	& \partial_M^\beta Q  = D \times \{0\} \cup 
	 \partial_M D\times [0,\infty) \cup  \{ \varpi \}.
\end{align}
In particular, $\partial_M^\beta Q$ does not depend on $\beta$.

\noindent
(ii)\ The $ \beta$-Martin compactification $Q^{\beta*}$ of $Q$ for
$ \partial_t+L$ is homeomorphic to the disjoint union of the 
topological product $ D^* \times [0,\infty)$ and the one point set 
$ \{ \varpi \}$, where a fundamental neighborhood system of $ \varpi$ is
given by the family 
$\{\varpi\} \cup D^* \times (N,\infty),\  N>1$.
In particular, $Q^{\beta*}$ does not depend on $\beta$.

\noindent
(iii)\ The $\beta$-Martin kernel $ K_\beta$ is given as follows:
For $(x,t) \in Q$,
\vphantom{$\displaystyle\int^A$}
\begin{align} \label{eq:4.9} 
	& K_\beta(x,t;y,0)=
	\frac{p(x,t; y,0)}{  \beta\left( p(\,\cdot\, ; y,0)\right)},
	\quad (y,0) \in D\times\{0\},\\
	& K_\beta(x,t;\xi,s)=
	\frac{q(x,t; \xi,s)}{  \beta\left( q(\,\cdot\, ; \xi,s)\right)},
	\quad (\xi,s) \in \partial_M D\times [0,\infty), \label{eq:4.10}
\end{align} 
and $K_\beta(x,t; \varpi)=0$.
\end{Prop}

This proposition is shown in the same way as Proposition 4.8 of
\cite{Mu07}.

\begin{Lem}\label{lem:4.6}
Let $(\xi,s) \in \left( \partial_M D \setminus \partial_m D \right)
\times [0,\infty)$. Then
there exists a finite Borel measure $\gamma$ on $\partial_M D$
supported by $\partial_m D$ such that
\begin{align} \label{eq:4.11} 
	& q(\,\cdot\, ; \xi,s)= \int_{\partial_m D}
	q(\,\cdot\, ; \eta,s)\, d \gamma(\eta).
\end{align}
\end{Lem}

\noindent{\bf Proof}\quad
For reader's convenience, we give a sketch of the proof for the case
$s=0$. (For details, see the proof of Lemma 4.10 of \cite{Mu07}.)
By the elliptic Martin representation theorem, 
there exists a unique finite 
Borel measure $\mu$ on $ \partial_M D$ supported by $ \partial_m D$ 
such that 
\[
	K(x,\xi)=\int_{\partial_m D} K(x,\eta)\, d \mu(\eta).
\]
This together with \rfe{eq:4.7} implies
\begin{align} \label{eq:4.12} 
	& \int_0^\infty q(x,\xi,t)\,dt =
	\int_{\partial_m D}\left(\int_0^\infty q(x,\eta,t)\,dt \right) d \gamma(\eta),
\end{align}
where $ d \gamma(\eta)= \left[ H(x^0,\xi)/H(x^0,\eta)\right]\,d \mu(\eta)$
with
\[
	H(x,\eta)=\int_0^\infty q(x,\eta,t)\,dt.
\]
For $ \alpha >0$, denote by $ G_\alpha$ the Green function of $ L+\alpha$
on $D$. By the resolvent equation and [SSP], we then have
\begin{align}  \label{eq:4.13}
	& \int_0^\infty 
	e^{- \alpha t} q(x,\eta,t)\,dt \\
	= & \int_0^\infty q(x,\eta,t)\,dt 
	-\alpha \int_D G_\alpha(x,z)
	\left(\int_0^\infty q(z,\eta,t)\,dt \right) \,m(z)d \nu(z),\nonumber
\end{align}
for any $\eta \in \partial_M D$. By combining \rfe{eq:4.12}
and \rfe{eq:4.13},
we get
\[
\int_0^\infty e^{- \alpha t} \left(\int_{\partial_m D}q(x,\eta,t)\,
d \gamma(\eta)\right)\,dt =
\int_0^\infty e^{- \alpha t} q(x,\xi,t)\,dt .
\]
Thus the Laplace transforms of $q(x,\xi,t)$
and $\int_{\partial_m D}q(x,\eta,t)\,d \gamma(\eta)$
coincide; and so \rfe{eq:4.11} holds. \qed

\begin{Lem}\label{lem:4.7}
Let $(\xi,s) \in \left( \partial_M D \setminus \partial_m D \right)
\times [0,\infty)$. Then $ q(\,\cdot\, ; \xi,s)$ is not minimal.
\end{Lem}

\noindent{\bf Proof}\quad
For reader's convenience, we give a proof. We have \rfe{eq:4.11}.
Suppose that $ q(\,\cdot\, ; \xi,s)$ is minimal.
Then, along the line given in the proof of Lemma 12.12 of
\cite{Hel69}, we obtain from \rfe{eq:4.11} that
the support of $\gamma$ consists of a single point.
Thus, for some $ \eta \in \partial_m D$ and constant $C$ 
\[
	q(\,\cdot\, ; \xi,s)=Cq(\,\cdot\, ; \eta,s).
\]
Hence, by Lemma~\ref{lem:4.4}, $ \xi=\eta$;
which is a contradiction. 
\qed

\begin{Lem}\label{lem:4.8}
Let $ (\xi,s) \in \partial_m D \times (0,\infty)$.
Then $q(\,\cdot\, ; \xi,s)$ is minimal if and only if
$q(\,\cdot\, ; \xi,0)$ is minimal.
\end{Lem}

\noindent{\bf Proof}\quad
Assume that $q(\,\cdot\, ; \xi,0)$ is minimal.
Suppose that a nonnegative solution $u$ of \rfe{eq:1.1}
satisfies $ u(\,\cdot\, ) \leq q(\,\cdot\, ; \xi,s)$
on $Q$. Put $ v(x,t)= u(x,t+s)$.
Then $ v(\,\cdot\, )\leq q(\,\cdot\, ; \xi,0)$.
Thus $ v(\,\cdot\, )= C q(\,\cdot\, ; \xi,0)$ for some constant
$C$. Hence $u(x,t) = C q(x,t; \xi,s)$ for $t>s$, and 
$u(x,t) = 0 =C q(x,t; \xi,s)$ for $t \le s$. This shows that
$q(\,\cdot\, ; \xi,s)$ is minimal.
Next, assume that $q(\,\cdot\, ; \xi,s)$ is minimal.
Suppose that a nonnegative solution $u$ of \rfe{eq:1.1}
satisfies $ u(\,\cdot\, ) \leq q(\,\cdot\, ; \xi,0)$
on $Q$. Put $ v(x,t)= u(x,t-s)$ for $t>s$, and
$v(x,t) = 0$ for $0<t \le s$. 
Then $v(\,\cdot\, )\leq q(\,\cdot\, ; \xi,s)$.
Thus $ v(\,\cdot\, )= C q(\,\cdot\, ; \xi,s)$ for some constant
$C$. Hence $u(x,t) = C q(x,t; \xi,0)$. This shows that
$q(\,\cdot\, ; \xi,0)$ is minimal.
\qed

By Theorem~\ref{th:4.1} and Lemmas~\ref{lem:4.7} and \ref{lem:4.8},
we have the following proposition.

\begin{Prop} \label{pro:4.9}
There exists a Borel subset $R$ of $\partial_M D$ such that
\[
	R \subset \partial_m D, \quad
	{\cal L}^\beta_m Q = R \times [0,\infty),
\]
for any $u \in P^0_\beta(Q)$ there exists a unique 
Borel measure $ \lambda$ on 
$ \partial_M D\times [0,\infty)$ 
which is supported by $R \times [0,\infty)$ and satisfies
\begin{align} \label{eq:4.14} 
	u(x,t) = \int_{R \times [0,\infty)} q(x,\xi,t-s)\, d \lambda(\xi,s)
	\quad (x,t) \in Q.
\end{align}
\end{Prop}

\begin{Lem}\label{lem:4.10}
Let $ (\xi,s) \in \partial_m D \times [0,\infty)$.
Then $q(\,\cdot\, ; \xi,s)$ is minimal.
\end{Lem}

\noindent{\bf Proof}\quad
Suppose that $q(\,\cdot\, ; \xi,0)$ is not minimal.
Then $\xi \notin R$ and 
\[
	q(x,\xi,t) =
	\int_{R \times [0,\infty)} q(x,\eta,t-s)\, d \lambda(\eta,s)
\]
for some Borel measure $ \lambda$. We have
\begin{align*} 
	K(x,\xi) \int_0^\infty q(x^0,\xi,t)\,dt 
	& = \int_0^\infty q(x,\xi,t) \,dt \\
	& = \int_{R \times [0,\infty)} d \lambda(\eta,s)
	K(x,\eta) \int_0^\infty q(x^0,\eta,t)\,dt .
\end{align*}
Thus
\[
K(x,\xi) = \int_{R} K(x,\eta) \,d \Lambda(\eta)
\]
for some Borel measure $\Lambda$. But
$\xi \in \partial_m D \setminus R \ $ and 
$\  R \subset \partial_m D$. This contradicts the uniqueness of a
representing measure in the elliptic Martin representation theorem.
Hence $q(\,\cdot\, ; \xi,0)$ is minimal; which together with 
Lemma~\ref{lem:4.8} shows Lemma~\ref{lem:4.10}.
\qed

\noindent{\bf Completion of the proof of Theorem~\ref{th:1.3}
in the case $I=(0,\infty)$}\quad
By Lemma~\ref{lem:4.10}, $R = \partial_m D$ and
\[
	 {\cal L}^\beta_m Q = \partial_m D \times [0,\infty).
\]
Thus Proposition~\ref{pro:4.9} shows Theorem~\ref{th:1.3}.
\qed

\section{Proof of Theorem~\ref{th:1.3}; the case $0<T<\infty$} \label{sec:5}

In this section we prove Theorem~\ref{th:1.3} in the case  
$0<T<\infty$ by making use of the results in Section 4. To this end, the 
following proposition plays a crucial role.

\begin{Prop} \label{pro:5.1}
Let $\xi \in \partial_M D$ and $0 \le s <r< \infty$. Then
\begin{align} \label{eq:5.1}
	\int_D p(x,y,t-r) q(y,r;\xi,s)d \mu(y)
	= q(x,t;\xi,s),
	\quad
	x \in D,\ t >r, 
\end{align}
where $d \mu(y) = m(y)\, d \nu(y)$
\end{Prop}

\noindent{\bf Proof}\quad
We first show \rfe{eq:5.1} for $\xi \in \partial_m D$.
Define $u(x,t)$ by
\begin{align} 
	u(x,t) &= q(x,t;\xi,s), &0<t \le r,  \nonumber\\
	u(x,t) &= \int_D p(x,y,t-r) q(y,r;\xi,s)d \mu(y), \ 
	&r<t< \infty.\label{eq:5.2} 
\end{align}
(We call $u$ the minimal extension of $q$ from $t=r$.) Then
we see that 
$u$ is a nonnegative solution of 
$(\partial_t + L)u=0$ in $D \times (0,\infty)$ such that  
$u(\,\cdot\, ) \leq q(\,\cdot\, ; \xi,s)$ on $D \times (0,\infty)$.
By Lemma~\ref{lem:4.10}, $u(\,\cdot\, ) = C q(\,\cdot\, ; \xi,s)$
for some constant $C$. But $u(x,t) = q(x,t;\xi,s)$ for $0<t \le r$.
Thus $C=1$, and so $u(\,\cdot\, ) = q(\,\cdot\, ; \xi,s)$. 

Next, let $\xi \notin \partial_m D$. By Lemma~\ref{lem:4.6},
there exists a finite Borel measure 
$\gamma$ on $\partial_M D$
supported by $\partial_m D$ such that
\begin{align} \label{eq:5.3} 
	& q(\,\cdot\, ; \xi,s)= \int_{\partial_m D}
	q(\,\cdot\, ; \eta,s)\, d \gamma(\eta).
\end{align}
Thus
\begin{align*}
	& \int_D p(x,y,t-r) q(y,r;\xi,s)d \mu(y) \\
	= & \int_{\partial_m D} \, d \gamma(\eta) 
	\int_D p(x,y,t-r) q(y,r;\eta,s)d \mu(y) \\
	= & \int_{\partial_m D} q(x,t; \eta,s)\, d \gamma(\eta) \\
	= & q(x,t;\xi,s).
\end{align*}
This proves \rfe{eq:5.1}.
\qed

\begin{Lem}\label{lem:5.2}
Let $\xi,\eta \in \partial_M D$, $ 0 \leq s,r < T$ and $ C >0$. If
\begin{align} \label{eq:5.4}
	& q(x,\xi,t-s)=Cq(x,\eta,t-r),\quad x \in D, \ 0<t<T,
\end{align}
then $\xi=\eta,\ s=r$ and $C=1$.
\end{Lem}

\noindent{\bf Proof}\quad
Choose $u$ such that $\max (r,s) <u < T$, and construct minimal extensions 
of both sides of \rfe{eq:5.4} from $t=u$. Then, by \rfe{eq:5.1} we have
\[
	q(x,\xi,t-s)=Cq(x,\eta,t-r),\quad x \in D, \ 0<t<\infty .
\]
By Lemma~\ref{lem:4.4}, this implies that 
$\xi=\eta,\ s=r$ and $C=1$.
\qed

Now, let $\beta$ be a measure on $Q = D \times (0,T)$
defined by 
\[
 d \beta(x,t)= a(x)b(t)m(x)\,d \nu(x) dt.
\]
Here $a(x)$ is a nonnegative continuous function on $D$ as
described in the beginning of Section 4, and $b(t)$ is 
a nonnegative continuous function on $\mathbf{R}$
such that $ b(t)>0$ on $ (T/2,T)$ and 
$b(t)=0$ on $ \mathbf{R} \setminus (T/2,T)$.
Let
$ K_\beta(\,\cdot\, ;\Xi)$, $ \partial^\beta_M Q $, $ \partial^\beta_m Q $, 
and $ Q^{\beta*}$ be 
the $\beta$-Martin kernel, $\beta$-Martin boundary, minimal
$\beta$-Martin boundary, and 
$\beta$-Martin compactification for $(Q,\ \partial_t + L)$
with $Q = D \times (0,T)$, respectively.
The following proposition is an analogue of Proposition~\ref{pro:4.5},
and is shown in the same way.

\begin{Prop}\label{pro:5.3}
(i)\ The $\beta$-Martin boundary $\partial_M^\beta Q$ of $Q$ for 
$ \partial_t +L$ is equal to the disjoint union of $D \times \{0\}$,
$ \partial_M D\times [0,T)$ and the one point set $ \{ \varpi \}$:
\begin{align} \label{eq:5.5} 
	& \partial_M^\beta Q  = D \times \{0\} \cup 
	 \partial_M D\times [0,T) \cup  \{ \varpi \}.
\end{align}
In particular, $\partial_M^\beta Q$ does not depend on $\beta$.

\noindent
(ii)\ The $ \beta$-Martin compactification $Q^{\beta*}$ of $Q$ for
$ \partial_t+L$ is homeomorphic to the disjoint union of the 
topological product $ D^* \times [0,T)$ and the one point set 
$ \{ \varpi \}$, where a fundamental neighborhood system of $ \varpi$ is
given by the family
$\{\varpi\} \cup D^* \times (T- \varepsilon,T),\  0<\varepsilon <T/2$.
In particular, $Q^{\beta*}$ does not depend on $\beta$.

\noindent
(iii)\ The $\beta$-Martin kernel $ K_\beta$ is given as follows:
For $(x,t) \in Q$,
\vphantom{$\displaystyle\int^A$}
\begin{align} \label{eq:5.6} 
	& K_\beta(x,t;y,0)=
	\frac{p(x,t; y,0)}{  \beta\left( p(\,\cdot\, ; y,0)\right)},
	\quad (y,0) \in D\times\{0\},\\
	& K_\beta(x,t;\xi,s)=
	\frac{q(x,t; \xi,s)}{  \beta\left( q(\,\cdot\, ; \xi,s)\right)},
	\quad (\xi,s) \in \partial_M D\times [0,T), \label{eq:5.7}
\end{align} 
and $K_\beta(x,t; \varpi)=0$.
\end{Prop}

\begin{Lem}\label{lem:5.4}
Let $(\xi,s) \in \left( \partial_M D \setminus \partial_m D \right)
\times [0,T)$. Then $ q(\,\cdot\, ; \xi,s)$ is not minimal.
\end{Lem}

\noindent{\bf Proof}\quad
Suppose that $ q(\,\cdot\, ; \xi,s)$ is minimal. Then we obtain
from \rfe{eq:5.3} that
\[
	q(x,\xi,t-s)=Cq(x,\eta,t-s),\quad x \in D, \ 0<t<T,
\]
for some $\eta \in \partial_m D$ and $ C >0$. By 
Lemma~\ref{lem:5.2}, this is a contradiction.
\qed

\begin{Lem}\label{lem:5.5}
Let $(\xi,s) \in \partial_m D \times [0,T)$. 
Then $ q(\,\cdot\, ; \xi,s)$ is minimal.
\end{Lem}

\noindent{\bf Proof}\quad
Let $u$ be a nonnegative solution of 
$(\partial_t + L)u=0$ in $Q$ such that  
$u(\,\cdot\, ) \leq q(\,\cdot\, ; \xi,s)$ in $Q$.
For $r \in (s,T)$, let $u_r$ be the minimal extension of $u$
from $t=r$. By Proposition~\ref{pro:5.1},
\[
	u_r(x,t) \le q(x,t;\xi,s),\quad x \in D, \ t>0.
\]
By Lemma~\ref{lem:4.10}, there exists a constant $C_r$ such
that $u_r(x,t) = C_r q(x,t;\xi,s)$ for $t>0$. But
 $u_r(x,t) = u(x,t)$ for $0<t<r$. Thus $C_r$ is independent of 
$r$; and so $u(\,\cdot\, ) = C q(\,\cdot\, ; \xi,s)$ in $Q$ for
some constant $C$.
\qed

\noindent{\bf Completion of the proof of Theorem~\ref{th:1.3}
in the case $0 < T < \infty $ }\quad
Put 
\[
	{\cal L}^\beta_m Q = \partial^\beta_m Q \setminus
	\left( D \times \{0\} \cup \{ \varpi \}\right).
\]
By Proposition~\ref{pro:5.3}, 
Lemmas~\ref{lem:5.4} and \ref{lem:5.5}, we get
\[
	 {\cal L}^\beta_m Q = \partial_m D \times [0,T).
\]
Thus, Theorem 2.1 of \cite{Mu07} which is an analogue of 
Theorem~\ref{th:4.1} completes the proof.
\qed

\section{Integral representations; the case $I=(-\infty,0)$} \label{sec:6}

In this section we prove Theorem~\ref{th:1.4}. 
We begin with the following proposition, which can be shown in the same 
way as in the proof of Theorem 1 of \cite{BD89} (see also \cite{Pi04}).

\begin{Prop} \label{pro:6.1}
Assume [SIU]. Then
\begin{equation} \label{eq:6.1}
\lim_{t \to \infty} \frac{e^{\lambda_0 t} p(x,y,t)}{\phi_0(x) \phi_0(y)} 
= 1 \quad \text{uniformly in} \   (x,y) \in K \times D
\end{equation}
for any compact subset $K$ of $D$.
\end{Prop}

In the rest of this section we assume [SSP].
We may and shall assume that $ a = 0 <\lambda_0$.
By Theorem~\ref{th:1.1}, we have the following corollary
of Proposition~\ref{pro:6.1}.

\begin{Cor} \label{cor:6.2}
Assume [SSP]. Then,
for any compact subset $K$ of $D$ and $N>1$,
\begin{equation*} 
\lim_{s \to -\infty} \frac{p(x,y,t-s)}{e^{\lambda_0 s} \phi_0(y)} 
= e^{-\lambda_0 t} \phi_0(x) 
\  \text{uniformly in} \  (x,y,t) \in K \times D \times (-N,0).
\end{equation*}
\end{Cor}

\begin{Lem}\label{lem:6.3}
The solution $e^{-\lambda_0 t} \phi_0(x)$ is minimal.
\end{Lem}

\noindent{\bf Proof}\quad
Suppose that $e^{-\lambda_0 t} \phi_0(x)$ is not minimal.
Then, in view of Corollary~\ref{cor:6.2}, the same argument as in
the proof of Theorem~\ref{th:1.3} shows that 
for any nonnegative solution $u$ of the equation
\begin{equation*} 
(\partial_t + L)u=0 \quad \text{in} \quad Q = D \times (-\infty,0) 
\end{equation*}
there exists a unique 
Borel measure $\lambda$ on $\partial_M D \times (-\infty,0)$
supported by the set  $\partial_m D \times (-\infty,0)$
such that
\begin{equation*} 
u(x,t) = 
\int_{\partial_M D \times (-\infty,t)} q(x,\xi,t-s) d\lambda(\xi,s),
\quad  (x,t) \in Q.
\end{equation*}
Thus
\begin{equation} \label{eq:6.2}
e^{-\lambda_0 t} \phi_0(x) = 
\int_{\partial_M D \times (-\infty,t)} q(x,\xi,t-s) d\lambda(\xi,s),
\quad  (x,t) \in Q ,
\end{equation}
for such a measure $\lambda$. Now, fix $x$. It follows from 
Theorems~\ref{th:1.1} and \ref{th:1.2} that for any 
$\delta > 0$ there exists a positive constant $C_{\delta}$
such that 
\begin{equation} \label{eq:6.3}
{C_{\delta}}^{-1} \le 
\frac{q(x,\xi,\tau)}{e^{-\lambda_0 \tau} \phi_0(x)}
\le C_{\delta},
\quad  \tau \ge \delta, \   \xi \in \partial_M D .
\end{equation}
By (\ref{eq:4.4}),
\begin{equation} \label{eq:6.4} 
q(x,\xi,\tau) \leq C e^{-\alpha/\tau} \phi_0(x),
\quad \xi \in \partial_M D, \ 0<\tau <1,
\end{equation}
for some positive constants $\alpha$ and $C$.
By (\ref{eq:6.2}) and (\ref{eq:6.3}),
\begin{equation*} 
e^{\lambda_0} \phi_0(x) \ge
\int_{\partial_M D \times (-\infty,-2)} 
C_1^{-1} e^{-\lambda_0 (-1-s)} d\lambda(\xi,s). 
\end{equation*}
Thus
\begin{equation} \label{eq:6.5} 
\int_{\partial_M D \times (-\infty,-2)} e^{\lambda_0 s} d\lambda(\xi,s)
\le C_1 \phi_0(x).
\end{equation}
For $t<-2$ and $0<\delta <1$, we have
\begin{equation} \label{eq:6.6}
\phi_0(x) = 
\int_{\partial_M D \times \{ (-\infty,t-\delta] \cup (t-\delta,t)\} } 
e^{\lambda_0 (t-s)} q(x,\xi,t-s) e^{\lambda_0 s} d\lambda(\xi,s).
\end{equation}
In view of (\ref{eq:6.4}) and (\ref{eq:6.5}), we choose $\delta$
so small that the integral on 
$\partial_M D \times (t-\delta,t)$ of the right hand side of 
(\ref{eq:6.6}) is smaller than $\phi_0(x)/3$. Then, in view of
(\ref{eq:6.3}) and (\ref{eq:6.5}), we choose $t<-2$ with
$|t|$ being so large that the integral on 
$\partial_M D \times (-\infty,t-\delta]$ of the right hand side of 
(\ref{eq:6.6}) is smaller than $\phi_0(x)/3$. This is a contradiction.
\qed

\noindent{\bf Completion of the proof of Theorem~\ref{th:1.4} }\quad
By virtue of Corollary~\ref{cor:6.2} and  
Lemma~\ref{lem:6.3}, the same argument as in the proof of 
Theorem~\ref{th:1.3} shows Theorem~\ref{th:1.4}.
\qed

\section{Examples} \label{sec:7}

In this section  we give two examples in order to illustrate a scope of 
Theorem~\ref{th:1.3}.  Throughout this section $L_0$ is a uniformly 
elliptic operator on $\mathbf{R}^n$ of the form 
\[
L_0 u= -\sum_{i,j=1}^n\,\partial_i\,\left(\,a_{ij}(x)\,\partial_{j} u\,\right),
\]
where $a(x)=\left[\, a_{ij}(x)\,\right]_{i,j=1}^n$ is a symmetric matrix-valued 
measurable function on $\mathbf{R}^n$ satisfying, for some $\Lambda >0$, 
\[\Lambda^{-1}\,|\xi|^2\,\leq\, \sum_{i,j=1}^n\,a_{ij}(x)\,\xi_{i}\xi_{j} 
\leq  \Lambda\,|\xi|^2 , \quad x, \xi \in  \mathbf{R}^n.\]  

\noindent {\bf 7.1. }  Let $V(x)$  be a measurable function in 
$L^{\infty}_{\mathrm{loc}}(\,\mathbf{R}^n\,)$, and $L=L_0+V(x)$ on $D=\mathbf{R}^n$.

\begin{Thm}\label{Tm:7.1} Suppose that there exist a positive constant $c<1$ 
and a positive continuous increasing function $\rho$ on $[0,\infty)$ such that
\begin{equation} 
c\,\left[\, \rho(|x|)\,\right]^2 \,\leq \, V(x) \,\leq \,  
\left[\, \rho(|x|)\,\right]^2, \quad x \in \mathbf{R}^n, \label{eq:7.1}
\end{equation}
\begin{equation} \label{eq:7.2}
c \,\rho\left(\,r+\frac{c}{\rho(r)}\,\right)\,\leq\,   \rho(r),\quad  r\geq 0.
\end{equation}
Assume that
\begin{equation} 
\int_{1}^{\infty} \frac{dr}{\rho(r)}< \infty.\label{eq:7.3}
\end{equation}
Then $1$ is a small perturbation of $L$ on $\mathbf{R}^n$. Thus Theorem~\ref{th:1.3} 
holds true.
\end{Thm}
{\bf Remark.} Compare this theorem with a non-uniqueness theorem of \cite{Mu94}.

\vskip 2mm
\noindent
{\bf Proof}\quad
We first note that \rfe{eq:7.2} yields
\[
c \rho(r) \leq 
c \rho\left(\,r - \frac{c}{\rho(r)} + \frac{c}{\rho\left(\,r- \frac{c}{\rho(r)}\,\right)}\,\right) 
\leq \rho\left(\,r-\frac{c}{\rho(r)}\,\right),
\qquad r \geq \frac{c}{\rho(0)} ,
\]
since $\rho$ is increasing.
We show the theorem by using the same approach as
in the proof of Theorem 5.1 of \cite{Mu02}. 
Put $b=c^{-2}$ and 
\[
\ell = \inf \{ j \in \mathbf{Z}; \rho(0) < b^j \}.
\]
For $k \geq \ell$,  \, put  \, 
$r_k = \sup \{ r \geq 0 ;  \, \rho(r) \leq b^k \}$.
By the continuity of $\rho$ and \rfe{eq:7.3}, $\rho(r_k) = b^k$ and 
 $\lim_{k \to \infty} r_k = \infty$.
By \rfe{eq:7.2},
\[
\rho(r_k + c b^{-k}) \leq c^{-1} \rho(r_k) = b^{1/2} b^k
< b^{k+1} = \rho(r_{k+1}).
\]
Thus \,$r_k + c b^{-k} < r_{k+1}$\, for $k \geq \ell$.
Define a positive continuously differentiabe increasing function
$\widetilde{\rho}$ on $[0,\infty)$ as follows:
Put $\widetilde{\rho}(r)= b^{\ell}$ for $r \leq r_{\ell}$, 
\[
\widetilde{\rho}(r)= b^{k+1} \quad
\text{for} \quad r_k + c b^{-k} \leq r \leq r_{k+1}
\quad (k \geq \ell);
\]
and 
$\widetilde{\rho}(r)= {\rho}_k (r)$ for 
$r_k \leq r \leq r_k + c b^{-k}$ \, 
$(k \geq \ell)$ \, by choosing a continuously differentiabe function 
${\rho}_k$ on
$[r_k,r_k + c b^{-k}]$ such that 
\[
{\rho}_k (r_k) = b^k, \quad {{\rho}_k}' (r_k) = 0, \quad
{\rho}_k (r_k + c b^{-k}) = b^{k+1}, \quad
{{\rho}_k}' (r_k + c b^{-k}) = 0,
\]
and 
\[
0 \leq {{\rho}_k}' (r) \leq B \, b^{2k}, \quad
r_k \leq r \leq r_k + c b^{-k},
\]
for some constant $B>0$ independent of $k$.
Then we have
\begin{equation}  \label{eq:7.4}
C^{-1}\leq \frac{\widetilde{\rho}(r)}{\rho(r)}\leq C, \quad
0 \leq {\widetilde{\rho}\, }'(r)
\leq C {\rho(r)}^2, \qquad r \geq 0,
\end{equation}
for some positive constant $C$.
Introduce a Riemannian metric $g=(g_{ij})_{i,j=1}^n$ by $g_{ij}
=\widetilde{\rho}(|x|)^2  \delta_{ij}$. Then $M=\mathbf{R}^n$ 
with this metric $g$ becomes a complete 
Riemannian manifold .  Furthermore, by \rfe{eq:7.2} and \rfe{eq:7.4}, $M$ has the 
bounded geometry property (1.1) of \cite{An97}.
The associated gradient $\nabla$ 
and divergence $\mathrm{div}$ are written as 
\[
\nabla= \widetilde{\rho}(|x|)^{-2}\,\nabla^0, \qquad
\mathrm{div}= \widetilde{\rho}(|x|)^{-n} \circ 
\mathrm{div}^0 \circ \widetilde{\rho}(|x|)^{n},
\]
where $\nabla^0$ and $\mathrm{div}^0$
are the standard gradient and divergence on $\mathbf{R}^n$. Put 
\begin{eqnarray*}
\cL&=&\widetilde{\rho}(|x|)^{-2}\,L, \\
m(x)=\widetilde{\rho}(|x|)^{2-n}, \quad
A(x)&=&\left[\,a_{ij}(x)\,\right]_{i,j=1}^n, \quad
\gamma(x)=\widetilde{\rho}(|x|)^{-2} \,V(x).
\end{eqnarray*}
Then 
\begin{eqnarray*}
{\cL}  u  &=&  - \frac{1}{m}\, \mathrm{div}\left(\,m A \nabla u \,\right)+ \gamma \\
&=& - \mathrm{div}\left( A \nabla u \,\right)-
\bigr{\langle} \frac{1}{m}\, A \, 
\nabla^0 m, \ \nabla u \bigr{\rangle}^0 + \gamma,
\end{eqnarray*}
where $\langle \cdot , \cdot {\rangle}^0$ is the standard inner product on $\mathbf{R}^n$.
Since the inner product $\langle \cdot , \cdot \rangle$ associated with the metric $g$ is 
written as 
\[
\langle X , Y \rangle  = \langle \widetilde{\rho}^{\,2} X , Y {\rangle}^0,
\]
we have
\begin{equation}  \label{eq:7.5}
{\cL}  u =
- \mathrm{div}\left( A \nabla u \,\right)-
\bigr{\langle} \widetilde{\rho}^{\, -2} 
\frac{ A \nabla^0 m}{m} , 
\ \nabla u \bigr{\rangle} + \gamma .
\end{equation}
By (\ref{eq:7.4}),
\[
|\nabla^0 m(x)|\,\leq\, C^3 \,|n-2|\,\widetilde{\rho}(|x|)\, m(x). 
\]
From this we have 
\begin{eqnarray*}
\bigr{\langle} \widetilde{\rho}^{\, -2} \frac{ A \nabla^0 m}{m} ,
\widetilde{\rho}^{\, -2} \frac{ A \nabla^0 m}{m} \bigr{\rangle} 
&\leq& \widetilde{\rho}^{\, -2} \Lambda^2 
(C^3 \,|n-2|\,\widetilde{\rho} \,)^2  \\
&\leq& \{  \Lambda (C^3\,|n-2|) \}^2 .
\end{eqnarray*}
By (\ref{eq:7.1}) and (\ref{eq:7.4}),
\[
c\,C^{-2}\,\leq\,\gamma(x)\,\leq\, C^2.
\]
Thus the operator ${\cL}-cC^{-2}/2$ has the Green function; and 
$\cL$ belongs to the class $\cD_{M}(\theta,\infty,\epsilon)$ introduced by 
Ancona \cite{An97}, where 
\[
\theta=\max\left(\,\Lambda, \Lambda (C^3 \,|n-2|) , C^2\,\right) ,
\quad  \epsilon=cC^{-2}/2. 
\]
Put 
\[
{\cL}_{2}   =  \widetilde{\rho}(|x|)^{-2}\,\left(\,L+1\,\right) 
= {\cL}+ \widetilde{\rho}(|x|)^{-2}.
\]
In order to apply the results of \cite{An97}, we proceed to estimate 
$\widetilde{\rho}(|x|)^{-2}$.  Let $d(x)$ be the Riemannian distance 
$\mathrm{dist}(0,x)$ from the origin $0$ to $x$, and put 
\[
\psi(r)=\int_0^r\, \widetilde{\rho}(s)\,ds.
\]
Then we see that $d(x)=\psi(|x|)$. Denote by $\psi^{-1}$ the inverse 
function of $\psi$, and put 
\[
\Phi(s)=\left[\,\widetilde{\rho}\left(\psi^{-1}(s)\,\right)\,\right]^{-2}, 
\quad s \geq 0.
\]
Then 
\[
0<\widetilde{\rho}(|x|)^{-2}=\Phi\left(\,d(x)\,\right), \quad x \in M.
\]
Furthermore, 
\begin{eqnarray*}
\int_{0}^\infty\,\Phi(s)\,ds &=& \int_{0}^\infty\,\Phi(\,\psi(r)\,)\,
\widetilde{\rho}(r)\,dr \\ 
&=&\int_{0}^\infty\,\frac{dr}{\widetilde{\rho}(r)}
\leq  C\,\int_{0}^\infty\,\frac{dr}{\rho(r)}\,dr\,<\infty.
\end{eqnarray*}
Hence, by virtue of Corollary 6.1, Theorems 1 and 2 of \cite{An97}, 
$\widetilde{\rho}(|x|)^{-2}$ is a small perturbation of $\cL$ on the manifold $M$. 
That is, 
for any $\varepsilon > 0$ there exists a compact subset $K$
of $D=M$ such that 
\begin{equation*}
\int_{D \setminus K} H(x,z) \widetilde{\rho}(|z|)^{-2}\,H(z,y) \, 
\widetilde{\rho}(|z|)^{n} dz
\le \varepsilon H(x,y), \qquad
x,y \in D \setminus K ,
\end{equation*}
where $dz$ is the Lebesgue measure on  $\mathbf{R}^n$, and 
$H(x,z)$ is the Green function of $\cL$ on $D$ with respect to the measure
$\widetilde{\rho}(|z|)^{n} dz$.  Denote by $G(x,z)$ the Green function of 
$L$ on $D$ with respect to the measure $dz$.
Since 
$\cL={\widetilde{\rho}(|x|)}^{-2} \, L $, 
we have
\[
H(x,z)= G(x,z)\, \widetilde{\rho}(|z|)^{2-n}
\]
Thus
\begin{equation*}
\int_{D \setminus K} G(x,z) \widetilde{\rho}(|z|)^{(2-n)-2}\,G(z,y) \, 
\widetilde{\rho}(|y|)^{2-n} \,  \widetilde{\rho}(|z|)^{n} dz
\leq \varepsilon G(x,y) \widetilde{\rho}(|y|)^{2-n} 
\end{equation*}
for any $x,y \in D \setminus K$.
Hence $1$ is a small perturbation of $L$ on $\mathbf{R}^n$. 
$\quad\square$

\vskip 2mm
\noindent
{\bf Remark.}  A sufficient 
condition for (\ref{eq:7.2}) is the following:
$\rho$ is a positive differentiable function on $[0, \infty)$ satisfying
\begin{equation} \label{eq:7.6}
0 \leq \rho'(r)\rho(r)^{-2} \leq C, \quad r \geq 0,
\end{equation}
for some positive constant $C$.  Indeed, from (\ref{eq:7.6}) we have 
\[
X(\delta) \equiv \rho\left(r+\frac{\delta}{\rho(r)}\right) \rho(r)^{-1} \leq
\exp [C \delta X(\delta)], \quad  r\geq 0, \  \delta >0.
\]
Put $\delta = (2Ce)^{-1}$, and let $\gamma \in (1,e)$ be the solution of
the equation 
\[
\exp [X/2e] = X.
\]
Then we get 
$1 \leq X(\delta) \leq \gamma$. Thus 
(\ref{eq:7.2}) holds with $c= \min(\delta, 1/\gamma)$.

\vskip 2mm
The condition (\ref{eq:7.3}) is sharp, since Theorem 6.2 of \cite{IM01} yields 
the following uniqueness theorem.

\begin{Thm}\label{Tm:7.2} Suppose that there exists  a positive continuous 
increasing function $\rho$ on $[0,\infty)$ such that
\begin{equation} 
 |V(x)| \leq   \rho(|x|)^2, \quad x \in \mathbf{R}^n. \label{eq:7.7}
\end{equation}
Assume that
\begin{equation} 
\int_{1}^{\infty} \frac{dr}{\rho(r)}= \infty.\label{eq:7.8}
\end{equation}
Then [UP] holds. Thus Fact AT holds true.
\end{Thm}

{\bf 7.2.} Throughout this subsection we assume that $D$ is a bounded domain of 
$\mathbf{R}^n$. Let $L$ be an elliptic operator on $D$ of the form 
\[L=\frac{1}{w(x)}\,L_0,\] 
where $w$ is a positive measurable function on $D$ such that 
$w, w^{-1} \in L^{\infty}_{\mathrm{loc}}( D )$.
\begin{Thm}\label{Tm:7.3} 
Let $D$ be a Lipschitz domain. Suppose that there 
exists a positive function $\psi$ on $(0,\infty)$ such that
$s^2\psi(s)$ is increasing and 
\begin{equation} \label{eq:7.9}
w(x)\leq \psi\left(\,\delta_D(x)\,\right), \quad x \in D, 
\end{equation}
where $\delta_D(x)=\, \mathrm{dist}\,(x,\partial D)$. Assume that
\begin{equation} \label{eq:7.10}
\int_{0}^{1}\, s\,\psi(s)\,ds\,<\,\infty.
\end{equation}
Then $1$ is a small perturbation of $L$ on $D$. Thus Theorem~\ref{th:1.3} 
holds true.
\end{Thm}

\noindent{\bf Remark.} (i) The first assertion of this theorem is implicitly shown in 
\cite{IM01} (see Theorem 7.11 and Remark 7.12 (ii) there).

(ii) The Lipschitz regularity of the domain $D$ is assumed only for 
the Hardy inequality to hold for any function in $C^{\infty}_{0}(D)$. Thus, 
for this theorem to hold, it suffices to assume (for example) that $D$ is 
uniformly $\Delta$-regular John domain or a simply connected domain of $\mathbf{R}^2$ 
(see  \cite{An86}, \cite{An97}).  

\vskip 2mm
\noindent
{\bf Proof of Theorem~\ref{Tm:7.3} } For $x \in D$, put 
\[D_x=\left\{\, y \in D; \,|x-y|<\frac{\delta_D(x)}{2}\right\}.\]
Then 
\[\frac{1}{2}\,\delta_D(x)\,\leq\,\delta_D(y)\,\leq\,\frac{3}{2}\,\delta_{D}(x),
\quad \, y \in D_x.\]
Thus
\begin{eqnarray*}
\delta_D(x)^2\,w(y) &\leq& 
4\,\delta_D(y)^2\,\psi\left(\,\delta_D(y)\,\right) \\
&\leq& 4\,\left(\,\frac{3}{2}\,\delta_{D}(x)\,\right)^2\,
\psi\left(\,\frac{3}{2}\,\delta_{D}(x)\,\right).
\end{eqnarray*}
Put $\Psi(s)=9 s^2 \psi\left((3/2)\,s\,\right)$. 
Then $\Psi(s)$ is increasing, and satisfies
\[
\delta_D(x)^2\,\left(\,\sup_{y \in D_x}\,w(y)\,\right)\leq 
\Psi\left(\,\delta_D(x)\,\right),
\qquad
 \int_{0}^{1} \frac{\Psi(s)}{s}\, ds\,< \infty.
\]
Hence, by virtue of Proposition 9.2, 
Theorem 9.1' and Corollary 6.1 of \cite{An97}, $w$ is a small 
perturbation of $L_0$ on $D$.
This implies that $1$ is a small 
perturbation of $L$ on $D$.
$\quad\square$

The condition (\ref{eq:7.10}) is sharp, since Theorem 7.8 and Lemma 7.6 of 
\cite{IM01} yield the following uniqueness theorem.

\begin{Thm}\label{Tm:7.4} Suppose that there exists  a positive continuous 
increasing function $\psi$ on $(0,\infty)$ such that
\begin{equation} 
c \psi\left(\, \delta_{D}(x)\,\right)\, \leq\, w(x) \,\leq\,  
\psi\left(\, \delta_{D}(x)\,\right) , \quad x \in D \label{eq:7.11}
\end{equation}
for some positive constant $c$, and 
\begin{equation} 
\nu \leq \frac{ \psi\left(\, \eta\,s\,\right)}{\psi(s)}\, \leq\, \nu^{-1} ,
\quad s>0, \, \, \frac{1}{2}\leq \eta\leq 2,  \label{eq:7.12}
\end{equation}
for some positive constant $\nu$. Assume 
\begin{equation} 
 \int_{0}^{1}\, \left[\,\psi(s)\,\left(\,\inf_{s\leq r \leq 1} r^2\,
\psi(r)\,\right)\,\right]^{\frac{1}{2}}\, ds \, = \, \infty. \label{eq:7.13}
\end{equation}
Then [UP] holds. Thus Fact AT holds true.
\end{Thm}


\end{document}